\documentclass[12pt, a4paper]{article}

\usepackage{a4wide}
\usepackage{amsmath}
\usepackage{amsfonts}
\usepackage{enumerate}
\date{}
\begin{document}
\title{\large\textbf{ON AN ABSTRACT FORMULATION OF A \\ THEOREM OF SIERPINSKI }}
\author{\normalsize
\textbf{D.SEN \& S.BASU}\\
}
\date{}	
\maketitle
{\small \noindent \textbf{{\textbf{AMS subject classification(2010):}}}} {\small $28A05$, $28A99$, $03E05$, $03E10$, $28D99$.\\}
{\small\noindent \textbf{\textbf{\textmd{\textbf{{Key words and phrases :}}}}} $G$-invariant class, $G$-invariant $k$-small system, $k$-additive measurable structure, admissible $k$-additive algebra, saturated set.}\\
\vspace{.03cm}\\
{\normalsize
\textbf{\textbf{ABSTRACT:}} In an earlier paper, we gave an abstract formulation of a theorem of Sierpi$\acute{n}$ski in uncountable commutative groups. In this paper, we prove a result which generalizes the earlier formulation.\\

\section{\large{INTRODUCTION}}
\normalsize Sierpi$\acute{n}$ski $[7]$ in one of his classical papers proved that there exist two Lebesgue measure zero sets in $\mathbb{R}$ whose algebraic sum is nonmeasurable. In establishing this result, he used Hamel basis and Steinhaus famous theorem on distance set. Several generalizations of Sierpi$\acute{n}$ski's theorem are available in the literature. Kharazishvili $[8]$ proved that for every $\sigma$-ideal $\mathcal I$ in $\mathbb{R}$ which is not closed with respect to the algebraic sum, and, for every $\sigma$-algebra $\mathcal S$($\supseteq\mathcal I$)for which the quotient algebra satisfies countable chain condition, there exist $X$, $Y$ $\in\mathcal I$ such that $X+Y\notin \mathcal S$. Now instead of the real line $\mathbb{R}$, if we choose a commutative group $G$ and any non-zero, $\sigma$-finite, complete, $G$-invariant (or, $G$-quasiinvariant)measure $\mu$, then an analogue of Sierpin$\acute{n}$ski's theorem can be established with respect to some extension of $\mu$. In fact, it was shown by Kharazishvili $[9]$ that for every uncountable commutative group $G$ and for any $\sigma$-finite, left $G$-invariant (or, $G$-quasiinvariant)measure $\mu$ on $G$, there exists a left $G$-invariant (or, $G$-quasiinvariant) complete measure $\mu^{\prime}$ extending $\mu$ and two sets $A,B\in \mathcal I(\mu^{\prime})$ (the $\sigma$-ideal of $\mu^{\prime}$-measure zero sets) such that $A+B\notin$ dom($\mu^{\prime}$). In earlier paper $[2]$, the present authors gave an abstract and generalized formulation of Sierpin$\acute{n}$ski's theorem in uncountable commutative groups which do not involve any use of measure.\
\vspace{.2cm}

Most of the notations, definitions and results of this paper are taken from $[2]$ (see also $[3]$, $[4]$). Throughout the paper, we identify every infinite cardinal with the least ordinal representing it, write card($E$) for the cardinality of any set $E$, and, use symbols such as $\xi$, $\rho$, $\alpha$, $k$ etc for any arbitrary infinite cardinal $k$ and $k^{^{+}}$ for the successor of $k$. Further, given an infinite group $G$ and a set $A\subseteq G$, we denote by $gA$ ($g\in G$) the set $\{gx : x\in A\}$ and call a class $\mathcal C$ of subsets of $G$ as $G$-invariant if $gA\in\mathcal C$ for every $g\in G$ and $A\in \mathcal C$.\\
\vspace{.2cm}

\textbf{DEFINITION $\textbf{1.1}$ :} A pair $(\Sigma,\mathcal I)$ consisting of two non-empty classes of subsets of $G$ is called a $G$-invariant, $k$-additive measurable structure on $G$ if \\
(i) $\Sigma$ is a algebra and $\mathcal I$ ($\subseteq\Sigma$) is a proper ideal in $G$.\\
(ii) Both $\Sigma$ and $\mathcal I$ are $k$-additive. This means that both the classes $\Sigma$ and $\mathcal I$ are closed with respect to the union of atmost $k$ number of sets.\\
(iii) $\Sigma$ and $\mathcal I$ are $G$-invariant.\
\vspace{.2cm}

A $k$-additive algebra $\Sigma$ is diffused if $\{x\}\in\Sigma$ for every $x\in G$ and a $k$-additive measurable structure $(\Sigma,\mathcal I)$ is called  $k^{^{+}}$-saturated if the cardinality of any arbitrary collection of mutually disjoint sets from $\Sigma\setminus\mathcal I$ is atmost $k$.\
\vspace{.2cm}

In the sixtees, Riecan and Neubrunn developed the notion of small systems and used the same to give abstract formulations of several well-known theorems in classical measure and integration (see $[13]$, $[14]$, $[15]$, etc) small system have been used by several other authors in the subsequent periods ($[5]$, $[6]$, $[11]$, $[12]$). The following Definition introduces a modified and generalized version of the same.\\
\vspace{.2cm}

\textbf{DEFINITION $\textbf{1.2}$ :} For any infinite cardinal $k$, a transfinite $k$-sequence $\{\mathcal N_{\alpha}\}_{_{\alpha<k}}$, of non empty classes of sets in $G$ is called a $G$-invariant, $k$-small system on $G$ if \\
(i) $\emptyset\in \mathcal N_{\alpha}$ for all $\alpha<k$.\\
(ii) Each $\mathcal N_{\alpha}$ is a $G$-invariant class.\\
(iii) $E\in \mathcal N_{\alpha}$ and $F\subseteq E$ implies $F\in \mathcal N_{\alpha}$\\
(iv) $E\in \mathcal N_{\alpha}$ and $F\in{\displaystyle{\bigcap_{\alpha<k}}\hspace{.01cm}{\mathcal N_{\alpha}}}$ implies $E\cup F\in\mathcal N_{\alpha}$\\
(v) For any $\alpha < k$, there exists $\alpha^{\ast}> \alpha$ such that for any one-to-one correspondence $\beta\rightarrow \mathcal N_{_{\beta}}$ with $\beta > \alpha^{\ast}$, ${\displaystyle{\bigcup_{\beta}{E_{_{\beta}}}}}\in \mathcal N_{_{\alpha}}$ whenever $E_{_{\beta}}\in \mathcal N_{_{\beta}}$.\\
(vi)  For any $\alpha , \beta < k$, there exists $\gamma > \alpha , \beta$  such that $\mathcal N_{_{\gamma}} \subseteq \mathcal N_{ _{\alpha}}$ and  $\mathcal N _{_{\gamma}} \subseteq \mathcal N _{_{\beta}}$.\
\vspace{.1cm}

We further define\\
\vspace{.2cm}

\textbf{DEFINITION $\textbf{1.3}$ :} A $G$-invariant $k$-additive algebra $\mathcal S$ on $G$ as admissible with respect to the $k$-small system $\{\mathcal N_{\alpha}\}{_{_{\alpha<k}}}$ if for every $\alpha<k$\\
(i) $\mathcal S\setminus\mathcal N_{\alpha}\neq\emptyset\neq\mathcal S\cap\mathcal N_{\alpha}$.\\
(ii) $\mathcal N_{\alpha}$ has a $S$-base i.e $E\in\mathcal N_{\alpha}$ is contained in some $F\in \mathcal N_{\alpha}\cap\mathcal S$,\\
and (iii) $\mathcal S\setminus\mathcal N_{\alpha}$ satisfies the $k$-chain condition, i.e, the cardinality of any arbitrary collection of mutually disjoint sets from $\mathcal S\setminus\mathcal N_{\alpha}$ is atmost $k$.\
\vspace{.1cm}

The above two Definitions have been used by the present authors in some of their recently done works (for example, see $[2]$, $[3]$, $[4]$).\ We set $\mathcal N_{\infty}= {\displaystyle{\bigcap_{\alpha<k}}}\hspace{.02cm}{\mathcal N_{\alpha}}$. From conditions (ii), (iii) and (v) of Definition $1.2$, it follows that $\mathcal N_{\infty}$ is a $G$-invariant, $k$-additive ideal in $G$ and denote by $\mathcal{\widetilde S}$ the $G$-invariant $k$-additive algebra generated by $\mathcal S$ and $\mathcal N_{\infty}$. Every element of $\mathcal{\widetilde S}$ is of the form $(X\setminus Y)\cup Z$ where $X\in\mathcal S$ and $Y,Z\in \mathcal N_{\infty}$ and $(\mathcal{\widetilde S}, \mathcal N_{\infty})$ turns out to be a $G$-invariant, $k$-additive measurable structure on $G$.\\
Moreover,\\
\vspace{.1cm}

\textbf{THEOREM $\textbf{1.4}$ :} \textit{If $\mathcal S$ is admissible with respect to $\{\mathcal N_{\alpha}\}{_{_{\alpha<k}}}$, then the $G$-invariant, $k$-additive measurable structure $(\mathcal{\widetilde S}, \mathcal N_{\infty})$ on G is $k^{^{+}}$-saturated.}\
\vspace{.2cm}

A proof of the above theorem follows directly from condition (iv) of Definition $1.2$ and conditions (i), (ii) and (iii) of Definition $1.3$, or in short from the admissibility of $\mathcal S$. Based on the above Definitions and theorems, some combinatorial properties of sets (Ch $7$, $[7]$) and an important representation theorem for infinite commutative groups (Appendix $2$, $[7]$), the present authors have proved in $[2]$ the following theorem.\\
\vspace{.2cm}

\textbf{THEOREM $\textbf{1.5}$ :} \textit{Let $G$ be an uncountable commutative group with card($G$)=$k^{^{+}}$. Let $\{\mathcal N_{\alpha}\}{_{_{\alpha<k}}}$ be a $G$-invariant, $k$-small system on $G$ and $\mathcal S$ be a diffused, $k$-additive algebra on $G$ which is also admissible with respect to $\{\mathcal N_{\alpha}\}{_{_{\alpha<k}}}$. Then there exists a subset $A$ of $G$ such that $A\in\mathcal N_{\infty}$ but $A+A\notin\mathcal{\widetilde S}$.
}
\section{\large{RESULT}}
\normalsize Theorem $1.5$ is an abstract formulation of Sierpinski's theorem given in terms of any diffused, $G$-invariant, $k$-additive measurable structure on a commutative group $G$ to which we have referred to in the introduction. In this section we prove a result which extends our previous formulation to groups that are not necessarily commutative.\\
\vspace{.1cm}

\textbf{DEFINITION $\textbf{2.1}[1]$ :} Let $\mathcal R$ be an equivalence relation on a set $X$ and $E\subseteq X$. The saturation of $E$ in $X$ with respect to the equivalence relation is the union of all equivalence classes of $\mathcal R$ whose intersection with $E$ is nonvoid.\\ In otherwords, it is  $\bigcup \{C  : C\cap E\neq \emptyset\hspace{.1cm} and\hspace{.1cm} C\in {X/\mathcal R}\}$\
\vspace{.1cm}

It is easy to check that if $H$ is a normal subgroup of any group $G$, then the saturation of any set $E$ in $G$ with respect to the equivalence relation generated by the quotient group $G/H$ is the set $HE$. If $E$ coincides with its saturation, then it is called saturated. Thus $E$ is saturated if $HE = E$. A saturated set is also called $H$-invariant $[10]$.\\
\vspace{.1cm}

\textbf{THEOREM $\textbf{2.2}$ :} \textit{Let $G$ be any uncountable group with card ($G$)=$k^{^{+}}$. Let $\{\mathcal N_{\alpha}\}{_{_{\alpha<k}}}$ be a $G$-invariant, $k$-small system on $G$ and $\mathcal S$ be a $G$-invariant, $k$-additive algebra on $G$ which is admissible with respect to $\{\mathcal N_{\alpha}\}{_{_{\alpha<k}}}$. We further assume that $G$ has a normal subgroup $H\in\mathcal S$ such that $G/H$ is commutative with card ($G/H$)=$k^{^{+}}$ and the saturation of any set $E$ in $G$ with respect to $G/H$ also belongs to $\mathcal S$.}\
\vspace{.01cm}

\textit{Then there exists a subset $A$ of $G$ such that $A\in\mathcal N_{\infty}$ and $AA\notin\mathcal{\widetilde S}$. }\\
\vspace{.2cm}

\textbf{PROOF : }We write $\Gamma = {G/H}$. By hypothesis $\Gamma$ is  commutative. Let $f:G\rightarrow \Gamma$ be the canonical homomorphism. We set $\mathcal S^{\prime}= \{Y\subseteq \Gamma : f^{-1}(Y)\in \mathcal S\}$ and $\mathcal N^{\prime}_{\alpha}= \{Y\subseteq\Gamma : f^{-1}(Y)\in\mathcal N_{\alpha}\}$ for any $\alpha<k$.\
\vspace{.1cm}

Since $\mathcal S$ is a $G$-invariant, $k$-additive algebra on $G$ and $f$ is a canonical homomorphism, so $\mathcal S^{\prime}$ is a $\Gamma$-invariant, $k$-additive algebra on $\Gamma$. Also since $H\in \mathcal S$, therefore $\mathcal S^{\prime}$ is diffused.\
\vspace{.1cm}

Condition (i) of Definition $1.2$ for $\{\mathcal N^{\prime}_{\alpha}\}{_{_{\alpha<k}}}$ is obvious. Let $h\in\Gamma$ and $F\in\mathcal N^{\prime}_{\alpha}$. Then $h=f(x)$ for every $x\in{gH}$ where $g\in G$ and $f^{-1}(F)\in \mathcal N_{\alpha}$. Since $\mathcal N_{\alpha}$ is $G$-invariant, therefore $f^{-1}(hF)= xf^{-1}(F)\in \mathcal N_{\alpha}$. Hence $hF\in \mathcal N^{\prime}_{\alpha}$ which proves condition (ii) of Definition $1.2$ for $\{\mathcal N^{\prime}_{\alpha}\}{_{_{\alpha<k}}}$. Finally, from the Definition of $\mathcal N^{\prime}_{\alpha}$ and some simple properties of inverse images of any function, it follows that conditions (iii)-(vi) of Definition $1.2$ also holds for $\{\mathcal N^{\prime}_{\alpha}\}{_{_{\alpha<k}}}$. Thus $\{\mathcal N^{\prime}_{\alpha}\}{_{_{\alpha<k}}}$ is a $\Gamma$-invariant, $k$-small system on $\Gamma$.\
\vspace{.1cm}

We shall now show that $\mathcal S^{\prime}$ is admissible with respect to $\{\mathcal N^{\prime}_{\alpha}\}{_{_{\alpha<k}}}$. Clearly, $\emptyset\in \mathcal S^{\prime}\cap\mathcal N^{\prime}_{\alpha}$ for $\alpha<k$. Since $\mathcal S$ is admissible with respect to $\{\mathcal N_{\alpha}\}{_{_{\alpha<k}}}$, so by (i) of Definition $1.3$, there exists for every $\alpha<k$, a set $A_{\alpha}\in\mathcal S\setminus\mathcal N_{\alpha}$. If $A_{\alpha}$ is saturated with respect to equivalence relation generated by the quotient group $G/H$, then $A_{\alpha}=f^{-1}(B_{\alpha})$ for some $B_{\alpha}\in{\mathcal S^{\prime}\setminus\mathcal N^{\prime}_{\alpha}}$. If $A_{\alpha}$ is not saturated, we replace it by $HA_{\alpha}$ which is saturated, and choose $B_{\alpha}$ such that $HA_{\alpha}=f^{-1}(B_{\alpha})$. Consequently $B_{\alpha}\in{\mathcal S^{\prime}\setminus\mathcal N^{\prime}_{\alpha}}$ and condition (i) of Definition $1.3$ is satisfied.\
\vspace{.1cm}

Let $F\in\mathcal N^{\prime}_{\alpha}$ and $E=f^{-1}(F)$. Then $E\in\mathcal N_{\alpha}$ by (ii) of Definition $1.3$ there exists $A\in\mathcal S\cap\mathcal N_{\alpha}$ such that $E\subseteq A$. If $A$ is saturated, then $A=f^{-1}(B)$ for some $B\in\mathcal S^{\prime}\cap\mathcal N^{\prime}_{\alpha}$ and $F\subseteq B$. If $A$ is not saturated, we choose the saturation of $G\setminus A$ i.e $H(G\setminus A)$ with respect to the equivalence relation generated by the quotient group $G/H$. But $H(G\setminus A)\in \mathcal S$ and so $G\setminus H(G\setminus A)\in \mathcal S$. Moreover, $G\setminus H(G\setminus A)$ is a subset of $A$. Therefore $G\setminus H(G\setminus A)\in\mathcal N_{\alpha}\cap\mathcal S$. We choose $B(\subseteq\Gamma)$ such that $G\setminus H(G\setminus A)=f^{-1}(B)$. Then $F\subseteq B$ and $B\in \mathcal S^{\prime}\cap\mathcal N^{\prime}_{\alpha}$. This shows that $\mathcal N^{\prime}_{\alpha}$ has a $\mathcal S^{\prime}$-base for every $\alpha<k$ and condition (ii) of Definition $1.3$ is proved. Lastly, any arbitrary collection of mutually disjoint sets from $\mathcal S^{\prime}\setminus\mathcal N^{\prime}_{\alpha}$ is atmost $k$ which follows directly from the fact that a similar result is true for the sets from $\mathcal S\setminus\mathcal N_{\alpha}$. This shows that $\mathcal S^{\prime}\setminus\mathcal N^{\prime}_{\alpha}$ satisfies the $k$-chain condition for every $\alpha<k$ which proves (iii) of Definition $1.3$.\
\vspace{.2cm}

Thus we find that $\mathcal S^{\prime}$ is a $\Gamma$-invariant, $k$-additive algebra on $\Gamma$ which is diffused and admissible with respect to the $\Gamma$-invariant, $k$-small system $\{\mathcal N^{\prime}_{\alpha}\}{_{_{\alpha<k}}}$ on $\Gamma$.\
\vspace{.2cm}

Let $\mathcal N^{\prime}_{\infty}= {\displaystyle{\bigcap_{\alpha<k}}}\hspace{.02cm}{\mathcal N^{\prime}_{\alpha}}$ and $\widetilde{\mathcal S^{\prime}}$ be the $\Gamma$-invariant, $k$-additive algebra generated by $\mathcal S^{\prime}$ and $\mathcal N^{\prime}_{\infty}$. Thus $(\widetilde{\mathcal S^{\prime}},\mathcal N^{\prime}_{\infty})$ is a $\Gamma$-invariant, $k$-additive, measurable structure on the quotient group $\Gamma$ which is $k^{^{+}}$-saturated. Hence by Theorem $1.5$, there exists $B\in\mathcal N^{\prime}_{\infty}$ such that $BB\notin\widetilde{\mathcal S^{\prime}}$. Let $A=f^{-1}(B)$. Then $AA=f^{-1}(B)f^{-1}(B)=f^{-1}(BB)$. So $AA$ is saturated. If possible, let $AA\in\widetilde{\mathcal S}$. Then $AA=E\Delta P$ where $E\in \mathcal S$, $P\in\mathcal N_{\infty}$ and $E$, $P$ are both saturated. Hence $E=f^{-1}(F)$, $P=f^{-1}(Q)$ where $F\in\mathcal S^{\prime}$, $Q\in\mathcal N^{\prime}_{\infty}$ and therefore $AA=E\Delta P= f^{-1}(F)\Delta f^{-1}(Q)=f^{-1}(F\Delta Q)=f^{-1}(BB)$. But this implies that $BB\in\widetilde{\mathcal S^{\prime}}$- a contradiction.\
\vspace{.2cm}

\textbf{REMARKS :} In general for Theorem $2.2$, $G$ need not be commutative. Let $H^{\prime}$ be a noncommutative group with card($H^{\prime}$)= $\omega$ (the first infinite cardinal) and $A^{\prime}$ be a commutative group with card($A^{\prime}$)= $\omega_{1}$ (the first uncountable cardinal). We set $G= H^{\prime}\times A^{\prime}$ as the external direct product of $H^{\prime}$ and $A^{\prime}$. Then $G$ is isomorphic with the internal direct product $HA$ where $H=\{(h,e_{A^{\prime}}): h\in H^{\prime}\}$ and $A=\{(e_{_{H^{\prime}}},a): a\in A^{\prime}\}$. Moreover $G$ is noncommutative having $H$ as a normal subgroup and $G/H$= $A$ is commutative with card($G/H$)= $\omega_{1}$.
\begin{center}

\end{center}

\textbf{\underline{AUTHOR'S ADDRESS}}\

{\normalsize
\textbf{S.Basu}\\
\vspace{.02cm}
\hspace{.28cm}
\textbf{Dept of Mathematics}\\
\vspace{.02cm}
\hspace{.35cm}
\textbf{Bethune College, Kolkata} \\
\vspace{.01cm}
\hspace{.3cm}
\textbf{W.B. India}\\
\vspace{.02cm}
\hspace{.1cm}
\hspace{.29cm}\textbf{{e-mail : sanjibbasu08@gmail.com}}\\
%\vspace{.0001cm}

\textbf{D.Sen}\\
\vspace{.02cm}
\hspace{.35cm}
\textbf{Saptagram Adarsha vidyapith (High), Habra, $\textbf{24}$ Parganas (North)} \\
\vspace{.01cm}
\hspace{.35cm}
\textbf{W.B. India}\\
\vspace{.02cm}
\hspace{.1cm}
\hspace{.29cm}\textbf{{e-mail : reachtodebasish@gmail.com}}
}

\end{document}